\documentclass{amsart}

\usepackage{xypic}
\input xy
\xyoption{all}
\usepackage{epsfig}
\usepackage{amsthm}
\usepackage{amssymb}
\usepackage{amsmath}
\usepackage{amscd}
\usepackage{color}
\usepackage[T1]{fontenc}

\newcommand{\bg}{\begin{equation}}
\newcommand{\ed}{\end{equation}}
\newcommand{\bga}{\begin{eqnarray}}
\newcommand{\eda}{\end{eqnarray}}
\newcommand{\pf}{\textbf{Proof:\ }}

\def\cbdu{\par{\raggedleft$\Box$\par}}

\newtheorem {Theorem}  {Theorem}

\numberwithin{Theorem}{section}

\newtheorem {Lemma}[Theorem]  {Lemma}

\theoremstyle{definition}
\newtheorem{Definition}[Theorem]{Definition}
\theoremstyle{remark}

\expandafter\chardef\csname pre amssym.def
at\endcsname=\the\catcode`\@ \catcode`\@=11
\def\undefine#1{\let#1\undefined}
\def\newsymbol#1#2#3#4#5{\let\next@\relax
 \ifnum#2=\@ne\let\next@\msafam@\else
 \ifnum#2=\tw@\let\next@\msbfam@\fi\fi
 \mathchardef#1="#3\next@#4#5}
\def\mathhexbox@#1#2#3{\relax
 \ifmmode\mathpalette{}{\m@th\mathchar"#1#2#3}%
 \else\leavevmode\hbox{$\m@th\mathchar"#1#2#3$}\fi}
\def\hexnumber@#1{\ifcase#1 0\or 1\or 2\or 3\or 4\or 5\or 6\or 7\or 8\or
 9\or A\or B\or C\or D\or E\or F\fi}

\font\teneufm=eufm10 \font\seveneufm=eufm7 \font\fiveeufm=eufm5
\newfam\eufmfam
\textfont\eufmfam=\teneufm \scriptfont\eufmfam=\seveneufm
\scriptscriptfont\eufmfam=\fiveeufm

\catcode`\@=\csname pre amssym.def at\endcsname

\newcounter{remark}
\newcommand{\R}{\mathbf{R}}

\def \ls {{\lambda_q^{2s}}}

\def  \R   {{\mathbb R}}

\def  \12  {{\frac{1}{2}}}

\def\build#1_#2^#3{\mathrel{\mathop{\kern 0pt#1}\limits_{#2}^{#3}}}

\begin{document}
%\currannalsline{0}{2006}

\title[Local well-posedness for Hall-MHD]{Local well-posedness for the Hall-MHD system in optimal Sobolev Spaces}

%\author{hello}

\author [Mimi Dai]{Mimi Dai}
\address{Department of Mathematics, Stat. and Comp.Sci., University of Illinois Chicago, Chicago, IL 60607,USA}
\email{mdai@uic.edu}

%%%use \Proof instead of \begin{proof}

%%%% use \Endproof instead of \end{proof}

%%%% use \references {999} instead of \begin{thebibliography}{99}

%%%%used \Endrefs instead of \end{thebibliography}

\begin{abstract}
We show that the viscous resistive magneto-hydrodynamics system with Hall effect is locally well-posed in $H^s(\R^n)\times H^{s+1-\varepsilon}(\R^n)$ with $s>\frac{n}2-1$ and any small enough $\varepsilon>0$ such that $s+1-\varepsilon>\frac{n}2$. This space is to date the largest local well-posedness space in the class of Sobolev spaces for the system. It is also optimal according to the predominant scalings of the two equations in the system.
%In this treatise we present the semigroup approach.... that has been developed over the last ten years... It emphasizes the dynamic viewpoint and is sufficiently general and flexible to encompass a great variety of concrete systems...  To date it is the only general method that applies to....

\bigskip


\hspace{0.02cm}CLASSIFICATION CODE: 76D03, 35Q35.
\end{abstract}

\maketitle

\section{Introduction}

Considered in this treatise is the three dimensional incompressible viscous resistive 
Hall-magneto-hydrodynamics (Hall-MHD) system:
\begin{equation}\label{HMHD}
\begin{split}
u_t+u\cdot\nabla u-b\cdot\nabla b+\nabla p-\nu\Delta u=0,\\
b_t+u\cdot\nabla b-b\cdot\nabla u+\eta\nabla\times((\nabla\times b)\times b)-\mu\Delta b=0,\\
\nabla \cdot u=0, 
\end{split}
\end{equation}
accompanied with the initial conditions
\begin{equation}\label{initial}
u(x,0)=u_0(x),\qquad b(x,0)=b_0(x), \qquad \nabla\cdot u_0=\nabla\cdot b_0=0,
\end{equation}
for $x\in\mathbb{R}^3$ and $t\geq 0$. In the system, $u$ represents the fluid velocity, $p$ is the fluid pressure and $b$
stands for the magnetic field. The parameters $\nu, \mu$ and $\eta$ denote the fluid viscosity, resistivity (electrical diffusivity) and the Hall effect coefficient, respectively.
%The parameter $\nu$ denotes the kinematic viscosity coefficient of the fluid and $\mu$ denotes the reciprocal of the magnetic Reynolds number. In this paper, we assume $\nu>0$ and $\mu>0$. 
It is important to observe that, if $\nabla\cdot b_0=0$, the divergence free condition for $b$ is propagated by the second equation of (\ref{HMHD}), see \cite{CL}.  
The Hall term $\nabla\times((\nabla\times b)\times b)$ distinguishes (\ref{HMHD}) from the usual MHD system (system (\ref{HMHD}) with $\eta=0$). 
In contrast to the latter one, the Hall-MHD model is more advantageous due to the fact that it can capture the essential characteristics of the magneto-hydrodynamics with strong magnetic reconnection where the Hall effect plays a significant role. Magnetic reconnection is a fundamental dynamical process in highly conductive plasmas in astrophysics, allowing for explosive and efficient magnetic to kinetic energy conversion.
For a more comprehensive physical background of the magnetic reconnection phenomena and the Hall-MHD model, we refer the readers to \cite{For, Li, SC} and references therein.

Despite its increasing popularity among the astrophysicists community,  the mathematical understanding of the Hall-MHD model is very limited. Conceptually, we can have a peek about the barriers from various perspectives. First, the Hall term launches new physics into the system at small length scales and hence intrinsically challenging into the mathematical analysis.
%Due to the inclusion of the Hall term, the model has an intrinsic multi-scale nature within which the interactions between small and large scales are delicate. 
Second, it is well-known that the main obstacle to understand the turbulent flows governed by the Navier-Stokes equation (NSE) relies on the nonlinearity such as $(u\cdot\nabla)u$. One can imagine that system (\ref{HMHD}) is more intricate than the NSE, for the former one contains the NSE and a magnetic field equation with the Hall term which appears more singular than  $(u\cdot\nabla)u$. Third, the natural scaling structure is a strong motivation in the study of both the NSE and the MHD system, who share the same scaling. However, the Hall term destroys such natural scaling. Into more details, for the MHD system, if $(u(x,t),p(x,t),b(x,t))$ solves (\ref{HMHD}) with $\eta=0$ with the initial data $(u_0(x),b_0(x))$, then the triplet $(u_\lambda(x,t),p_\lambda(x,t),b_\lambda(x,t))$ defined by
\begin{equation}\label{scale1}
u_\lambda(x,t)=\lambda u(\lambda x,\lambda^2t), \ \ p_\lambda(x,t)=\lambda^2 p(\lambda x,\lambda^2t), \ \ b_\lambda(x,t)=\lambda b(\lambda x,\lambda^2t)
\end{equation}
solves the same system with the data
\[u_{0\lambda}(x,t)=\lambda u_0(\lambda x), \ \ b_{0\lambda}(x,t)=\lambda b_0(\lambda x).\]
The scaling (\ref{scale1}) no longer holds for system (\ref{HMHD}) with $\eta>0$. On the other hand, we can extract the ``Hall equation'' 
\begin{equation}\notag%\label{Hall}
b_t+\nabla\times((\nabla\times b)\times b)=\Delta b
\end{equation}
which has the scaling 
\begin{equation}\label{scale-h}
b_\lambda(x,t)=b(\lambda x,\lambda^2t).
\end{equation}
Since the Hall term is the most singular nonlinearity in system (\ref{HMHD}), it suggests that the predominant scaling for (\ref{HMHD}) could be 
\begin{equation}\label{scale2}
u_\lambda(x,t)=\lambda u(\lambda x,\lambda^2t), \ \ p_\lambda(x,t)=\lambda^2 p(\lambda x,\lambda^2t), \ \ b_\lambda(x,t)= b(\lambda x,\lambda^2t).
\end{equation}
In fact, based on scaling (\ref{scale2}), we obtained a regularity criterion for (\ref{HMHD}) in three dimension which improves various criteria in the literature, see \cite{Dai-hmhd-reg}.

In this paper our interest is to find the largest possible (optimal) Sobolev space where system (\ref{HMHD}) is locally well-posed. On this topic, it was first shown in \cite{CWW} that system (\ref{HMHD}) in three dimension is locally well-posed in $H^s(\R^3)\times H^s(\R^3)$ with $s>\frac 52$. By taking (\ref{scale-h}) as the dominant scaling, in \cite{Dai-hmhd-well}, we obtained the local well-posedness of (\ref{HMHD}) in $H^s(\R^n)\times H^s(\R^n)$ with $s>\frac n2$. Even though the result of \cite{Dai-hmhd-well} improves that of \cite{CWW}, it seems that there is still room to have improvement, for the reason that the NSE is known to be locally well-posed in $H^s(\R^n)$ with $s>\frac n2-1$. In fact, motivated by scaling (\ref{scale2}), one expects that system (\ref{HMHD}) may be locally well-posed in $H^s(\R^n)\times H^{s+1}(\R^n)$ with $s>\frac n2-1$. In order to justify the conjecture, we need to treat the energy estimates for $u$ and $b$ separately, namely, $u$ in $H^s$ and $b$ in $H^{s+1}$. In this situation, we encounter the difficulty that no cancelation can be employed to deal with the two terms $b\cdot \nabla b$ and $b\cdot\nabla u$. To overcome this barrier, it comes to our mind that we need to optimize the estimates of the flux contributed from the two terms by fully employing the diffusion of both the $u$ and the $b$. Techniques based on the paradifferential calculus enables us to operate such optimizations. Surprisingly, it turns out that the local well-posedness space we can obtain is slightly larger than the conjectured one.
In deed, we prove the main result below.

\begin{Theorem}\label{thm}
Let $(u_0,b_0)\in H^s(\R^n)\times H^{s+1-\varepsilon}(\R^n)$ with $s>\frac{n}2-1$ and any small enough $\varepsilon>0$ such that $s+1-\varepsilon>\frac{n}2$. Assume $\nabla \cdot u_0=\nabla\cdot b_0=0$. There exists a time $T=T(\nu,\mu,\|u_0\|_{H^s},\|b_0\|_{H^{s+1-\varepsilon}})>0$ and a unique solution $(u,b)$ of (\ref{HMHD}) on $[0,T]$ such that 
\[(u,b)\in C([0,T); H^s(\R^n))\times C([0,T); H^{s+1-\varepsilon}(\R^n)).\]
\end{Theorem}

Regarding the result, the fact that $b$ needs to be in a space with higher regularity is determined by the Hall term. Predicted by the scaling (\ref{scale-h}) of the ``Hall equation'', the optimal Sobolev space of well-posedness for $b$ would be $H^{s+1}(\R^n)$ with $s>\frac{n}2-1$. However, as stated in Theorem \ref{thm}, the obtained well-posdness space for $b$ is $H^{s+1-\varepsilon}(\R^n)$ for any small $\varepsilon>0$. It may be explained by getting a closer look at the term $b\cdot\nabla u$. While 
estimating $\|b\cdot\nabla u\|_{H^r}$ by applying both diffusions of $u$ and $b$, it happens that we need to take $r$ slightly smaller than $s+1$.
%In fact, the Hall-MHD system and the MHD system satisfy the same basic energy law, which is a result of the cancelation in the flux from the Hall term. Such cancelation is due to special geometry arrangement of magnetic lines. It is different from the cancelation property in the convection terms ($u\cdot\nabla u$ and $u\cdot\nabla b$) which is due to the divergence free condition. 

%The rest of the paper is organized as follows: in Section \ref{sec:pre} we introduce some notations, recall the Littlewood-Paley decomposition theory briefly, and establish some auxiliary estimates to handle the Hall term; Section \ref{sec:reg} is devoted to proving Theorem \ref{thm}.  

\bigskip

\section{Preliminaries}
\label{sec:pre}

\subsection{Notation}
\label{sec:notation}
In order to avoid confusion, we specify a few notations. 
We denote by $A\lesssim B$ an estimate of the form $A\leq C B$ with
some absolute constant $C$, and by $A\sim B$ an estimate of the form $C_1
B\leq A\leq C_2 B$ with absolute constants $C_1$, $C_2$. For simplification, it is understood that
 $\|\cdot\|_p=\|\cdot\|_{L^p}$.
 % and $(\cdot, \cdot)$ for the $L^2$-inner product.

\subsection{Littlewood-Paley decomposition}
\label{sec:LPD}
As in our previous articles on the local well-posedness of magneto-hydrodynamics systems, the main tool is paradifferential calculus. 
%The techniques presented in this paper rely strongly on the Littlewood-Paley decomposition. 
%Thus we need certain amount of preparation by 
To be self-contained, we 
recall the Littlewood-Paley decomposition theory briefly, even though it appears in our earlier work on related topics. For a more detailed description on this theory we refer the readers to \cite{BCD} and \cite{Gr}. 

Let $\mathcal F$ and $\mathcal F^{-1}$ denote the Fourier transform and inverse Fourier transform, respectively. Define $\lambda_q=2^q$ for integers $q$. A nonnegative radial function $\chi\in C_0^\infty(\R^n)$ is chosen such that 
\begin{equation}\notag
\chi(\xi)=
\begin{cases}
1, \ \ \mbox { for } |\xi|\leq\frac{3}{4}\\
0, \ \ \mbox { for } |\xi|\geq 1.
\end{cases}
\end{equation}
Let 
\bg\notag
\varphi(\xi)=\chi(\frac{\xi}{2})-\chi(\xi)
\ed
and
\begin{equation}\notag
\varphi_q(\xi)=
\begin{cases}
\varphi(\lambda_q^{-1}\xi)  \ \ \ \mbox { for } q\geq 0,\\
\chi(\xi) \ \ \ \mbox { for } q=-1.
\end{cases}
\end{equation}
For a tempered distribution vector field $u$ we define the Littlewood-Paley projection
\begin{equation}\notag%\label{eq:LPP}
\begin{split}
&h=\mathcal F^{-1}\varphi, \qquad \tilde h=\mathcal F^{-1}\chi,\\
&u_q:=\Delta_qu=\mathcal F^{-1}(\varphi(\lambda_q^{-1}\xi)\mathcal Fu)=\lambda_q^n\int h(\lambda_qy)u(x-y)dy,  \qquad \mbox { for }  q\geq 0,\\
& u_{-1}=\mathcal F^{-1}(\chi(\xi)\mathcal Fu)=\int \tilde h(y)u(x-y)dy.
\end{split}
\end{equation}
By the Littlewood-Paley theory, the identity
\bg\notag
u=\sum_{q=-1}^\infty u_q
\ed
holds in the distributional sense. 
%Essentially the sequence of the smooth functions $\varphi_q$ forms a dyadic partition of the unit. 
For brevity, we agree with the notations
\bg\notag
u_{\leq Q}=\sum_{q=-1}^Qu_q, \qquad u_{(Q, N]}=\sum_{p=Q+1}^N u_p, \qquad \tilde u_q=\sum_{|p-q|\leq 1}u_p.
\ed
%Notice that by the definition of $\varphi_q$ 
%$$\supp (\hat u_q)\cap\supp (\hat u_p)=\emptyset  \qquad \mbox { if }  |p-q|\geq 2.$$

%By the Littlewood-Paley decomposition theory,
%we define the Besov spaces $B_{p,r}^{s}$.
%\begin{Definition}
%Let $s\in \mathbb R$, and $1\leq p, r\leq \infty$. Define the norm
%$$
%\|u\|_{ B_{p,r}^{s}}=\left(\sum_{q=-1}^\infty\lambda_q^{sr}\|u_q\|_p^r\right)^{\frac1r}.
%$$
%The Besov space $B_{p,r}^{s}$ is the space of tempered distributions $u$ such that the norm $\|u\|_{B_{p,r}^{s}}$ is finite.
%\end{Definition}
%Specially, 
\begin{Definition}
A tempered distribution $u$ belongs to the Besov space $ B_{p, \infty}^{s}$ if and only if
$$
\|u\|_{ B_{p, \infty}^{s}}=\sup_{q\geq-1}\lambda_q^s\|u_q\|_p<\infty.
$$
\end{Definition}
We can identify the Sobolev space $H^s$ by the Besov space $B^s_{2,2}$, i.e.
\[
  \|u\|_{H^s} \sim \left(\sum_{q=-1}^\infty\lambda_q^{2s}\|u_q\|_2^2\right)^{1/2}
\]
for each $u \in H^s$ and $s\in\R$.

%We recall Bernstein's inequality for the dyadic blocks of the Littlewood-Paley decomposition in the following.
\begin{Lemma}\label{le:bern} (Bernstein's inequality. See \cite{L}.)
Let $n$ be the space dimension and $r\geq s\geq 1$. Then for all tempered distributions $u$, we have
\bg\label{Bern}
\|u_q\|_{r}\lesssim \lambda_q^{n(\frac{1}{s}-\frac{1}{r})}\|u_q\|_{s}.
\ed
\end{Lemma}

\bigskip

%\subsection{Definition of solutions}
%\label{sec:sol}
%We recall some classical definitions of weak and regular solutions.
%\begin{Definition}\label{def:weak}
%A weak solution of (\ref{HMHD}) on $[0,T]$ (or $[0, \infty)$ if $T=\infty$) is a pair of functions $(u, b)$ in the class 
%$$
%u, b \in C_w([0,T]; L^2(\mathbb R^3))\cap L^2(0,T; H^1(\mathbb R^3)), %\quad b\in C_w([0,T]; L^2(\mathbb R^3)) \cap L^2(0,T; H^1(\mathbb R^3))
%$$
%with $u(0)=u_0, b(0)=b_0$, satisfying (\ref{HMHD}) in the distribution sense;
% and for all test functions $\phi\in C_0^\infty([0,T]\times\mathbb R^3)$ with $\nabla_x\cdot \phi=0$ {\color{red} Check the Hall term.}
%\begin{equation}\notag%\label{eq:uweak}
%\begin{split}
%&(u(t), \phi(t))-(u_0, \phi(0))\\
%=&\int_0^t(u(s), \partial_s\phi(s))+\nu(u(s), \Delta\phi(s))+ (u(s)\cdot\nabla\phi(s), u(s))\, ds,
%\end{split}
%\end{equation}
%\begin{equation}\notag%\label{eq:bweak}
%\begin{split}
%&(b(t), \phi(t))-(b_0, \phi(0))\\
%=&\int_0^t(b(s), \partial_s\phi(s))+\mu(b(s), \Delta\phi(s))+ (b(s)\cdot\nabla\phi(s), b(s))\, ds\\
%&+\int_0^t(\nabla\times b(s)\times b(s), \nabla\times \phi(s))\, ds,
%\end{split}
%\end{equation}
%moreover, the following energy inequality 
%\begin{equation}\notag%\label{energy}
%\begin{split}
%&\|u(t)\|_2^2+\|b(t)\|_2^2+2\nu\int_{t_0}^t\|\nabla u(s)\|_2^2ds
%+2\mu\int_{t_0}^t\|\nabla b(s)\|_2^2ds\\
%\leq &\|u(t_0)\|_2^2+\|b(t_0)\|_2^2
%\end{split}
%\end{equation}
%is satisfied for almost all $t_0\in(0, T)$ and all $t\in(t_0, T]$.
%\end{Definition}

\subsection{Bony's paraproduct and commutator}
\label{sec-para}

Bony's paraproduct formula %It is usually sufficient to use a weaker form of this formula:
\begin{equation}\label{Bony}
\begin{split}
\Delta_q(u\cdot\nabla v)=&\sum_{|q-p|\leq 2}\Delta_q(u_{\leq{p-2}}\cdot\nabla v_p)+
\sum_{|q-p|\leq 2}\Delta_q(u_{p}\cdot\nabla v_{\leq{p-2}})\\
&+\sum_{p\geq q-2} \Delta_q(\tilde u_p \cdot\nabla v_p),
\end{split}
\end{equation}
will be used constantly to decompose the nonlinear terms in energy estimate. 
We will also use the notation of the commutator
\begin{equation} \label{commudef}
[\Delta_q, u_{\leq{p-2}}\cdot\nabla]v_p:=\Delta_q(u_{\leq{p-2}}\cdot\nabla v_p)-u_{\leq{p-2}}\cdot\nabla \Delta_qv_p.
\end{equation}
\begin{Lemma}\label{le-commu}
The commutator satisfies the following estimate, for any $1<r<\infty$
\[\|[\Delta_q,u_{\leq{p-2}}\cdot\nabla] v_p\|_{r}\lesssim \|\nabla u_{\leq p-2}\|_\infty\|v_p\|_{r}.\]
\end{Lemma}

\subsection{Auxiliary estimates}
\label{sec-commu}

To handle the Hall term $\nabla \times((\nabla \times b)\times b)$, more preparation is needed. We first introduce two more commutators and their estimates.
We define that, for vector valued functions $F$ and $G$,
\begin{equation}\label{comm-v}
[\Delta_q,F\times\nabla\times]G=\Delta_q(F\times(\nabla\times G))-F\times(\nabla\times G_q),
\end{equation}
\begin{equation}\label{comm-v2}
[\Delta_q,\nabla\times F\times]G=\Delta_q(\nabla\times F\times G)-\nabla\times F\times G_q.
\end{equation}
In principle, the commutators will be used to reveal certain cancellation; and to shift derivative from high modes to low modes. It was shown in \cite{Dai-hmhd-reg} they satisfy the following estimates. 
\begin{Lemma}\label{le-Hall1}
Let $F$ and $G$ be vector valued functions. Assume $\nabla\cdot F=0$ and $F$, $G$ vanish at large $|x|\in \R^3$. For any $1<r<\infty$, we have
\[\|[\Delta_q,F\times\nabla\times]G\|_r\lesssim  \|\nabla F\|_\infty\|G\|_r;\]
\[\|[\Delta_q,\nabla\times F\times]G\|_r\lesssim  \|\nabla F\|_\infty\|G\|_r.\]
\end{Lemma}

\begin{Lemma}\label{le-Hall2}
Let $F$, $G$ and $H$ be vector valued functions. Assume $F$, $G$ and $H$ vanish at large $|x|\in \R^3$. For any $1<r_1, r_2<\infty$ with $\frac1{r_1}+\frac1{r_2}=1$, we have
\[\left|\int_{\R^3}[\Delta_q,\nabla\times F\times]G\cdot\nabla\times H\, dx\right|
\lesssim \|\nabla^2 F\|_\infty\|G\|_{r_1}\|H\|_{r_2}.\]
\end{Lemma}

\bigskip

\section{A priori estimate}
\label{sec-est}
In this section, we establish a priori estimate for smooth solutions in $H^s(\R^n)\times H^{r}(\R^n)$ with appropriate index $s$ and $r$. Such estimate is the most crucial ingredient in the argument of local well-posedness, which is rather standard for dissipative equations, see \cite{MB}. Thus we only present the following theorem and its proof. 

\begin{Theorem}\label{thm-priori} Let $(u_0,b_0)\in H^s(\R^n)\times H^{r}(\R^n)$ with $s>\frac n2-1$ and $\frac n2<r\leq s+1-\varepsilon$ for small enough $\varepsilon>0$. There exists a time $T=T(\nu,\mu,\|u_0\|_{H^s}, \|b_0\|_{H^r})>0$ such that the Hall-MHD system (\ref{HMHD}) has a solution $(u,b)$ satisfying
\[u\in L^\infty(0,T; H^s(\R^n))\cap L^2(0,T; H^{s+1}(\R^n)),\]
\[b\in L^\infty(0,T; H^r(\R^n))\cap L^2(0,T; H^{r+1}(\R^n)).\]
\end{Theorem}
The proof involves certain amount of computations and estimates which will be divided into several lemmas, each carrying an estimate for a flux term. To start,
multiplying the first equation of (\ref{HMHD}) by $\lambda_q^{2s}\Delta_qu_{q}$ and the second one by $\lambda_q^{2r}\Delta_qb_{q}$, and adding up for all $q\geq -1$,  we obtain 
\begin{equation}\label{ineq-uq}
\frac12\frac{d}{dt}\sum_{q\geq -1}\lambda_q^{2s}\|u_q\|_2^2
 +\nu\sum_{q\geq-1}\lambda_q^{2s+2}\|u_q\|_2^2\leq -I_1-I_2,
\end{equation}
\begin{equation}\label{ineq-bq}
\frac12\frac{d}{dt}\sum_{q\geq -1}\lambda_q^{2r}\|b_q\|_2^2 +\mu\sum_{q\geq-1}\lambda_q^{2r+2}\|b_q\|_2^2\leq -I_3-I_4-I_5,
\end{equation}
with
\begin{equation}\notag
\begin{split}
I_1=&\sum_{q\geq -1}\lambda_q^{2s}\int_{\R^3}\Delta_q(u\cdot\nabla u)\cdot u_q\, dx, \qquad
I_2=-\sum_{q\geq -1}\lambda_q^{2s}\int_{\R^3}\Delta_q(b\cdot\nabla b)\cdot u_q\, dx,\\
I_3=&\sum_{q\geq -1}\lambda_q^{2r}\int_{\R^3}\Delta_q(u\cdot\nabla b)\cdot b_q\, dx,\qquad
I_4=-\sum_{q\geq -1}\lambda_q^{2r}\int_{\R^3}\Delta_q(b\cdot\nabla u)\cdot b_q\, dx,\\
I_5=&-\sum_{q\geq -1}\lambda_q^{2r}\int_{\R^3}\Delta_q((\nabla\times b)\times b)\cdot \nabla\times b_q\, dx.
\end{split}
\end{equation}
To fully exploit cancelations in the flux terms $I_1$, $I_3$ and $I_5$, we will apply commutator estimates along with Bony's paraproduct and some fundamental inequalities. While $r\neq s$, there is no cancelation in $I_2+I_4$, and hence $I_2$ and $I_4$ will be treated in slightly different ways.

\begin{Lemma}\label{le-i1}
Let $s>\frac n2-1$. We have that, for some absolute constants $\gamma_1, \gamma_2>0$,
\begin{equation}\notag%\label{est-i1}
|I_1|\leq \frac{\nu}{8}\sum_{q\geq -1}\lambda_q^{2s+2}\|u_q\|_2^2+C_\nu\|u\|_{H^s}^{2+\gamma_1}
+C_\nu\|u\|_{H^s}^{2+\gamma_2}.
\end{equation}
\end{Lemma}
\pf
Using Bony's paraproduct (\ref{Bony}) followed by the commutator notation (\ref{commudef}), $I_1$ is decomposed as
\begin{equation}\notag%\label{eq:i2}
\begin{split}
I_1=
&\sum_{q\geq -1}\sum_{|q-p|\leq 2}\lambda_q^{2s}\int_{\R^3}\Delta_q(u_{\leq p-2}\cdot\nabla u_p)u_q\, dx\\
&+\sum_{q\geq -1}\sum_{|q-p|\leq 2}\lambda_q^{2s}\int_{\R^3}\Delta_q(u_{p}\cdot\nabla u_{\leq{p-2}})u_q\, dx\\
&+\sum_{q\geq -1}\sum_{p\geq q-2}\lambda_q^{2s}\int_{\R^3}\Delta_q(u_p\cdot\nabla\tilde u_p)u_q\, dx\\
=&I_{11}+I_{12}+I_{13},
\end{split}
\end{equation}
with 
\begin{equation}\notag
\begin{split}
I_{11}=&\sum_{q\geq -1}\sum_{|q-p|\leq 2}\lambda_q^{2s}\int_{\R^3}[\Delta_q, u_{\leq{p-2}}\cdot\nabla] u_pu_q\, dx\\
&+\sum_{q\geq -1}\sum_{|q-p|\leq 2}\lambda_q^{2s}\int_{\R^3}u_{\leq{q-2}}\cdot\nabla \Delta_q u_p u_q\, dx\\
&+\sum_{q\geq -1}\sum_{|q-p|\leq 2}\lambda_q^{2s}\int_{\R^3}(u_{\leq{p-2}}-u_{\leq{q-2}})\cdot\nabla\Delta_qu_p u_q\, dx\\
=&I_{111}+I_{112}+I_{113}.
\end{split}
\end{equation}
Thanks to the facts $\sum_{q-2\leq p\leq q+2}\Delta_qu_p=u_q$ and $\nabla\cdot u_{\leq q-2}=0$, the term $I_{112}$ vanishes. Notice that $I_{12}$ and $I_{13}$ can be treated in the analogous way as $I_{111}$ and $I_{113}$, respectively. Thus we will only show the estimates of $I_{111}$ and $I_{113}$.
Applying the commutator estimate in Lemma \ref{le-commu} and Bernstein's inequality to $I_{111}$ gives rise to
\begin{equation}\notag
\begin{split}
|I_{111}|\leq&\sum_{q\geq -1}\sum_{|p-q|\leq 2}\lambda_q^{2s}\|\nabla u_{\leq p-2}\|_\infty\|u_p\|_2\|u_q\|_2\\
\lesssim &\sum_{q\geq -1}\lambda_q^{2s}\|u_q\|_2^2\sum_{p\leq q}\lambda_p^{\frac n2+1}\|u_p\|_2\\
\lesssim &\sum_{q\geq -1}\lambda_q^{(s+1)\theta}\|u_q\|_2^\theta \lambda_q^{s(2-\theta)}\|u_q\|_2^{2-\theta}\sum_{p\leq q}\lambda_p^{(s+1)\delta}\|u_p\|_2^\delta\lambda_p^{s(1-\delta)}\|u_p\|_2^{1-\delta}\left(\lambda_q^{-\theta}\lambda_p^{\frac n2+1-s-\delta}\right)\\
\lesssim &\sum_{q\geq -1}\lambda_q^{(s+1)\theta}\|u_q\|_2^\theta \lambda_q^{s(2-\theta)}\|u_q\|_2^{2-\theta}\sum_{p\leq q}\lambda_p^{(s+1)\delta}\|u_p\|_2^\delta\lambda_p^{s(1-\delta)}\|u_p\|_2^{1-\delta}\lambda_{p-q}^\theta
\end{split}
\end{equation}
with constants $\theta$ and $\delta$ satisfying $0<\theta<2$, $0<\delta<1$ and 
\begin{equation}\label{para1}
s\geq \frac n2+1-\theta-\delta.
\end{equation}
It then follows from Young's inequality with $(r_1,r_2,r_3,r_4)\in (1,\infty)^4$ satisfying 
\begin{equation}\label{para2}
\frac1{r_1}+\frac1{r_2}+\frac1{r_3}+\frac1{r_4}=1, \ \ r_1=\frac2\theta, \ \ r_3=\frac2{\delta}
\end{equation}
such that for some $\theta_1>0, \theta_2>0$
\begin{equation}\notag
\begin{split}
|I_{111}|\leq &\frac{\nu}{64}\sum_{q\geq -1}\lambda_q^{2s+2}\|u_q\|_2^2+ C_{\nu}\sum_{q\geq -1}\left(\lambda_q^{2s}\|u_q\|_2^2\right)^{\frac{(2-\theta)r_2}{2}}\\
&+\frac{\nu}{64}\sum_{q\geq -1}\sum_{p\leq q}\lambda_p^{2s+2}\|u_p\|_2^2\lambda_{p-q}^{\theta_1}+C_\nu\sum_{q\geq -1}\sum_{p\leq q}\left(\lambda_p^{2s}\|u_p\|_2^2\right)^{\frac{(1-\delta)r_4}{2}}\lambda_{p-q}^{\theta_2}\\
\leq &\frac{\nu}{32}\sum_{q\geq -1}\lambda_q^{2s+2}\|u_q\|_2^2+C_\nu\left(\sum_{q\geq-1}\lambda_q^{2s}\|u_q\|_2^2\right)^{\frac{(2-\theta)r_2}{2}}
+C_\nu\left(\sum_{q\geq-1}\lambda_q^{2s}\|u_q\|_2^2\right)^{\frac{(1-\delta)r_4}{2}}
\end{split}
\end{equation}
Notice that (\ref{para1}) and (\ref{para2}) imply that $s>\frac n2-1$.

To estimate $I_{113}$, it follows from H\"older, Bernstein and Young's inequalities that
\begin{equation}\notag
\begin{split}
|I_{113}|\leq&\sum_{q\geq -1}\sum_{|p-q|\leq 2}\lambda_q^{2s}\|u_{\leq p-2}-u_{\leq q-2}\|_2\|\nabla u_p\|_\infty\|u_q\|_2\\
\lesssim &\sum_{q\geq -1}\lambda_q^{2s+\frac n2+1}\|u_q\|_2^3\\
\lesssim &\sum_{q\geq -1}\lambda_q^{(s+1)\theta}\|u_q\|_2^\theta \lambda_q^{s(3-\theta)}\|u_q\|_2^{3-\theta}\lambda_q^{\frac n2+1-s-\theta}\\
\leq &\frac \nu{32}\sum_{q\geq -1}\lambda_q^{2s+2}\|u_q\|_2^2+C_\nu\left(\sum_{q\geq-1}\lambda_q^{2s}\|u_q\|_2^2\right)^{\frac{3-\theta}{2-\theta}}
\end{split}
\end{equation}
for $s\geq \frac n2+1-\theta$ and $0<\theta<2$.  Thus 
\begin{equation}\notag%\label{est-i1}
I_1\leq \frac{\nu}{8}\|\nabla u\|_{H^s}^2+C_\nu \| u\|_{H^s}^{2+\gamma_1}+C_\nu \| u\|_{H^s}^{2+\gamma_2}
%|I_1|\leq \frac{\nu}{8}\sum_{q\geq -1}\lambda_q^{2s+2}\|u_q\|_2^2+C_\nu\left(\sum_{q\geq-1}\lambda_q^{2s}\|u_q\|_2^2\right)^{1+\gamma_1}+C_\nu\left(\sum_{q\geq-1}\lambda_q^{2s}\|u_q\|_2^2\right)^{1+\gamma_2}
\end{equation}
for $s>\frac n2-1$ and some $\gamma_1, \gamma_2>0$.

\cbdu

\begin{Lemma}\label{le-i2}
Let $\frac n2+s-2r\leq 0$ and $s<r$. The following estimate holds 
\begin{equation}\notag%\label{est-i2}
|I_2|\leq \frac{\nu}{8} \sum_{q\geq -1}\lambda_q^{2s+2}\|u_q\|_2^2+C_\nu\| b\|_{H^r}^{4}.
\end{equation}
%Notice that if $r=s+3-2\alpha$, then $s\geq \frac n2-6+4\alpha >\frac n2-1$ provided $\frac12<\alpha\leq \frac 54$.
\end{Lemma}
%{\color{red} Can we raise the power of $\|b\|_{H^r}$ from 4 to higher number? If so, we may obtain global regularity for small initial data.} 
\pf
We first decompose $I_2$ by using Bony's paraproduct, 
\begin{equation}\notag%\label{eq:i2}
\begin{split}
I_2=
&-\sum_{q\geq -1}\sum_{|q-p|\leq 2}\lambda_q^{2s}\int_{\R^3}\Delta_q(b_{\leq p-2}\cdot\nabla b_p)u_q\, dx\\
&-\sum_{q\geq -1}\sum_{|q-p|\leq 2}\lambda_q^{2s}\int_{\R^3}\Delta_q(b_{p}\cdot\nabla b_{\leq{p-2}})u_q\, dx\\
&-\sum_{q\geq -1}\sum_{p\geq q-2}\lambda_q^{2s}\int_{\R^3}\Delta_q(b_p\cdot\nabla\tilde b_p)u_q\, dx\\
=&I_{21}+I_{22}+I_{23}.
\end{split}
\end{equation}
%with 
%\begin{equation}\notag
%\begin{split}
%I_{21}=&-\sum_{q\geq -1}\sum_{|q-p|\leq 2}\lambda_q^{2s}\int_{\R^3}[\Delta_q, b_{\leq{p-2}}\cdot\nabla] b_pu_q\, dx\\
%&-\sum_{q\geq -1}\sum_{|q-p|\leq 2}\lambda_q^{2s}\int_{\R^3}b_{\leq{q-2}}\cdot\nabla \Delta_q b_p u_q\, dx\\
%&-\sum_{q\geq -1}\sum_{|q-p|\leq 2}\lambda_q^{2s}\int_{\R^3}(b_{\leq{p-2}}-b_{\leq{q-2}})\cdot\nabla\Delta_qb_p u_q\, dx\\
%=&I_{211}+I_{212}+I_{213}.
%\end{split}
%\end{equation}
Due to the lack of cancelation, $I_{21}$ is the worst term which can be estimated as
\begin{equation}\notag
\begin{split}
|I_{21}|\leq & \sum_{q\geq -1}\sum_{|q-p|\leq 2}\lambda_q^{2s+1}\|b_{\leq p-2}\|_\infty\|b_p\|_2\|u_q\|_2\\
\lesssim & \sum_{q\geq -1}\lambda_q^{2s+1}\|b_q\|_2\|u_q\|_2\sum_{p\leq q}\lambda_p^{\frac n2}\|b_p\|_2\\
\lesssim & \sum_{q\geq -1}\lambda_q^{s+1}\|u_q\|_2 \lambda_q^r\|b_q\|_2\sum_{p\leq q}\lambda_p^{r}\|b_p\|_2\lambda_{q-p}^{s-r}\lambda_p^{\frac n2+s-2r}\\
\lesssim & \sum_{q\geq -1}\lambda_q^{s+1}\|u_q\|_2 \lambda_q^r\|b_q\|_2\sum_{p\leq q}\lambda_p^{r}\|b_p\|_2\lambda_{q-p}^{s-r}
\end{split}
\end{equation}
for $\frac n2+s-2r\leq 0$.
As a result, Young's inequality gives rise to
\begin{equation}\notag
\begin{split}
|I_{21}|\leq &\frac{\nu}{16} \sum_{q\geq -1}\lambda_q^{2s+2}\|u_q\|_2^2+C_\nu\sum_{q\geq -1}\left(\lambda_q^r\|b_q\|_2\sum_{p\leq q}\lambda_p^{r}\|b_p\|_2\lambda_{q-p}^{s-r}\right)^2.
\end{split}
\end{equation}
Then we apply Jensen's inequality, if $s<r$, %and followed by changing the order of the summations
\begin{equation}\notag
\begin{split}
|I_{21}|\leq &\frac{\nu}{16} \sum_{q\geq -1}\lambda_q^{2s+2}\|u_q\|_2^2+C_\nu\sum_{q\geq -1}\lambda_q^{2r}\|b_q\|_2^2\sum_{p\leq q}\lambda_p^{2r}\|b_p\|_2^2\lambda_{q-p}^{s-r}\\
\leq &\frac{\nu}{16} \sum_{q\geq -1}\lambda_q^{2s+2}\|u_q\|_2^2+C_\nu\left(\sum_{ q\geq-1}\lambda_q^{2r}\|b_q\|_2^2\right)^{2}.
\end{split}
\end{equation}
We claim that $I_{22}$ shares the same estimate as $I_{21}$. Indeed, the following inequality holds
\begin{equation}\notag
\begin{split}
|I_{22}|\lesssim & \sum_{q\geq -1}\lambda_q^{2s}\|b_q\|_2\|u_q\|_2\sum_{p\leq q}\lambda_p^{\frac n2+1}\|b_p\|_2\lesssim |I_{21}|.
\end{split}
\end{equation}
To move the derivative from high modes to low modes in $I_{23}$, we apply integration by parts 
\begin{equation}\notag
\begin{split}
|I_{23}|=&\left|\sum_{q\geq -1}\sum_{p\geq q-2}\lambda_q^{2s}\int_{\R^3}\Delta_q(b_p\otimes\tilde b_p)\cdot\nabla u_q\, dx\right|.
\end{split}
\end{equation}
It then follows from H\"older's and Bernstein's inequalities 
\begin{equation}\notag
\begin{split}
|I_{23}|\lesssim & \sum_{q\geq -1}\lambda_q^{2s+1}\|u_q\|_2\sum_{p\geq q-4}\|b_p\|_2\|b_p\|_\infty\\
\lesssim & \sum_{q\geq -1}\lambda_q^{2s+1}\|u_q\|_2\sum_{p\geq q-4}\lambda_p^{\frac n2}\|b_p\|_2^2\\
\lesssim & \sum_{q\geq -1}\lambda_q^{s+1}\|u_q\|_2 \sum_{p\geq q-4}\lambda_p^{2r}\|b_p\|_2^2\lambda_{q-p}^{s}\lambda_p^{\frac n2+s-2r}\\
\lesssim & \sum_{q\geq -1}\lambda_q^{s+1}\|u_q\|_2 \sum_{p\geq q-4}\lambda_p^{2r}\|b_p\|_2^2\lambda_{q-p}^{s}
\end{split}
\end{equation}
for $\frac n2+s-2r\leq 0$. Applying Young's inequality, Jensen's inequality and changing order of the summations yields
\begin{equation}\notag
\begin{split}
|I_{23}|\leq &\frac{\nu}{16} \sum_{q\geq -1}\lambda_q^{2s+2}\|u_q\|_2^2+C_\nu\sum_{q\geq -1}\left(\sum_{p\geq q-4}\lambda_p^{2r}\|b_p\|_2^2\lambda_{q-p}^{s}\right)^2\\
\leq &\frac{\nu}{16} \sum_{q\geq -1}\lambda_q^{2s+2}\|u_q\|_2^2+C_\nu\sum_{q\geq -1}\sum_{p\geq q-4}\lambda_p^{4r}\|b_p\|_2^4\lambda_{q-p}^{s}\\
\leq &\frac{\nu}{16} \sum_{q\geq -1}\lambda_q^{2s+2}\|u_q\|_2^2+C_\nu\sum_{p\geq -1}\lambda_p^{4r}\|b_p\|_2^4\sum_{q\leq p+4}\lambda_{q-p}^{s}\\
\leq &\frac{\nu}{16} \sum_{q\geq -1}\lambda_q^{2s+2}\|u_q\|_2^2+C_\nu\left(\sum_{ q\geq-1}\lambda_q^{2r}\|b_q\|_2^2\right)^{2}.
\end{split}
\end{equation}
It completes the proof.
%To conclude, we obtain 
%\begin{equation}\notag%\label{est-i2}
%|I_2|\leq \frac{3\nu}{16}\|\nabla u\|_{H^s}^2+C_\nu \| b\|_{\dot H^r}^{4}+C_\nu \| b\|_{\dot H^s}^{2+\gamma}
%\frac{3\nu}{16} \sum_{q\geq -1}\lambda_q^{2s+2}\|u_q\|_2^2+C_\nu\left(\sum_{ q\geq-1}\lambda_q^{2r}\|b_q\|_2^2\right)^{2}+C_\nu\left(\sum_{ q\geq-1}\lambda_q^{2r}\|b_q\|_2^2\right)^{1+\gamma}
%\end{equation}
%for $\frac n2+s-2r\leq 0$ and $s<r$.

\cbdu

\begin{Lemma}\label{le-i3}
Let $s>\frac n2-1$ and $\frac n4+\frac s2<r<s+2-\varepsilon$ with small enough $\varepsilon>0$. We have the estimate  
\begin{equation}\notag%\label{est-i3-1}
\begin{split}
|I_{3}|
\leq &\frac{\nu}{8}\sum_{q\geq -1}\lambda_q^{2s+2}\|u_q\|_2^2+\frac{\mu}{8}\sum_{q\geq -1}\lambda_q^{2r+2}\|b_q\|_2^2
+C_{\nu,\mu}\| u\|_{H^s}^{2+\gamma_3}
+C_{\nu,\mu}\| b\|_{H^r}^{2+\gamma_4}
\end{split}
\end{equation}
for some constants $\gamma_3, \gamma_4>0$. 
\end{Lemma}
\pf
As for $I_1$, we first decompose $I_3$ by Bony's paraproduct 
\begin{equation}\notag%\label{eq:i2}
\begin{split}
I_3=
&\sum_{q\geq -1}\sum_{|q-p|\leq 2}\lambda_q^{2r}\int_{\R^3}\Delta_q(u_{\leq p-2}\cdot\nabla b_p)b_q\, dx\\
&+\sum_{q\geq -1}\sum_{|q-p|\leq 2}\lambda_q^{2r}\int_{\R^3}\Delta_q(u_{p}\cdot\nabla b_{\leq{p-2}})b_q\, dx\\
&+\sum_{q\geq -1}\sum_{p\geq q-2}\lambda_q^{2r}\int_{\R^3}\Delta_q(u_p\cdot\nabla\tilde b_p)b_q\, dx\\
=&I_{31}+I_{32}+I_{33},
\end{split}
\end{equation}
and further decompose $I_{31}$ by using the commutator to
\begin{equation}\notag
\begin{split}
I_{31}=&-\sum_{q\geq -1}\sum_{|q-p|\leq 2}\lambda_q^{2r}\int_{\R^3}[\Delta_q, u_{\leq{p-2}}\cdot\nabla] b_pb_q\, dx\\
&-\sum_{q\geq -1}\sum_{|q-p|\leq 2}\lambda_q^{2r}\int_{\R^3}u_{\leq{q-2}}\cdot\nabla \Delta_q b_p b_q\, dx\\
&-\sum_{q\geq -1}\sum_{|q-p|\leq 2}\lambda_q^{2r}\int_{\R^3}(u_{\leq{p-2}}-u_{\leq{q-2}})\cdot\nabla\Delta_qb_p b_q\, dx\\
=&I_{311}+I_{312}+I_{313}.
\end{split}
\end{equation}
It is not hard to see that $I_{312}=0$. % due to the fact $\sum_{|p-q|\leq 2}\Delta_qb_p=b_q$ and $\nabla\cdot u_{\leq q-2}=0$.
By the commutator estimate in Lemma \ref{le-commu}, we infer
\begin{equation}\notag
\begin{split}
|I_{311}|\leq&\sum_{q\geq -1}\sum_{|p-q|\leq 2}\lambda_q^{2r}\|\nabla u_{\leq p-2}\|_\infty\|b_p\|_2\|b_q\|_2\\
\lesssim &\sum_{q\geq -1}\lambda_q^{2r}\|b_q\|_2^2\sum_{p\leq q}\lambda_p^{\frac n2+1}\|u_p\|_2\\
\lesssim &\sum_{q\geq -1}\lambda_q^{(r+1)\theta}\|b_q\|_2^\theta \lambda_q^{r(2-\theta)}\|b_q\|_2^{2-\theta}\sum_{p\leq q}\lambda_p^{(s+1)\delta}\|u_p\|_2^\delta\lambda_p^{s(1-\delta)}\|u_p\|_2^{1-\delta}\left(\lambda_q^{-\theta}\lambda_p^{\frac n2+1-s-\delta}\right)\\
\lesssim &\sum_{q\geq -1}\lambda_q^{(r+1)\theta}\|b_q\|_2^\theta \lambda_q^{r(2-\theta)}\|b_q\|_2^{2-\theta}\sum_{p\leq q}\lambda_p^{(s+1)\delta}\|u_p\|_2^\delta\lambda_p^{s(1-\delta)}\|u_p\|_2^{1-\delta}\lambda_{p-q}^{\theta}\\
\end{split}
\end{equation}
for parameters $\theta$ and $\delta$ satisfying $0<\theta<2$, $0<\delta<1$ and 
\begin{equation}\label{para5}
s\geq \frac n2+1-\theta-\delta.
\end{equation}
It then follows from Young's inequality with $(r_1,r_2,r_3,r_4)\in (1,\infty)^4$ satisfying 
\begin{equation}\label{para6}
\frac1{r_1}+\frac1{r_2}+\frac1{r_3}+\frac1{r_4}=1, \ \ r_1=\frac2\theta, \ \ r_3=\frac2{\delta}
\end{equation}
such that 
\begin{equation}\notag
\begin{split}
|I_{311}|
%\leq &\frac{\nu}{64}\sum_{q\geq -1}\lambda_q^{2s+2}\|u_q\|_2^2+ C_{\nu}\sum_{q\geq -1}\left(\lambda_q^{2s}\|u_q\|_2^2\right)^{\frac{(2-\theta)r_2}{2}}\\&+\frac{\nu}{64}\sum_{q\geq -1}\sum_{p\leq q}\lambda_p^{2s+2}\|u_p\|_2^2\lambda_{p-q}^{\theta_1}+C_\nu\sum_{q\geq -1}\sum_{p\leq q}\left(\lambda_p^{2s}\|u_p\|_2^2\right)^{\frac{(1-\delta)r_4}{2}}\lambda_{p-q}^{\theta_2}\\
\leq &\frac{\nu}{32}\sum_{q\geq -1}\lambda_q^{2s+2}\|u_q\|_2^2+\frac{\mu}{32}\sum_{q\geq -1}\lambda_q^{2r+2}\|b_q\|_2^2\\
&+C_{\nu,\mu}\left(\sum_{q\geq-1}\lambda_q^{2s}\|u_q\|_2^2\right)^{1+\gamma_3}
+C_{\nu,\mu}\left(\sum_{q\geq-1}\lambda_q^{2r}\|b_q\|_2^2\right)^{1+\gamma_4}
\end{split}
\end{equation}
for some constants $\gamma_3,\gamma_4>0$.
Notice that (\ref{para5}) and (\ref{para6}) imply for large enough $r_2$ and $r_4$, and $\delta, \theta$ close enough to 1, there exists a small $\varepsilon>0$ such that
 \begin{equation}\notag%\label{s2}
 s\geq\frac n2-\theta+\varepsilon>\frac n2-1.
 \end{equation}
We observe that $|I_{313}|\lesssim |I_{311}|$, and hence $I_{313}$ enjoys the same estimate of $I_{311}$.
%\begin{equation}\notag
%\begin{split}
%|I_{313}|\leq&\sum_{q\geq -1}\sum_{|p-q|\leq 2}\lambda_q^{2r}\|u_{\leq p-2}-u_{\leq q-2}\|_2\|\nabla b_p\|_\infty\|b_q\|_2\\
%\lesssim &\sum_{q\geq -1}\lambda_q^{2r+\frac n2+1}\|u_q\|_2\|b_q\|_2^2\\
%\lesssim &\sum_{q\geq -1}\lambda_q^{s+1}\|u_q\|_2 \lambda_q^{2r}\|b_q\|_2^{2}\lambda_q^{\frac n2-s}\\
%\leq &\frac \nu{32}\sum_{q\geq -1}\lambda_q^{2s+2}\|u_q\|_2^2+C_\nu\left(\sum_{q\geq-1}\lambda_q^{2r}\|b_q\|_2^2\right)^2
%\end{split}
%\end{equation}
%for $\frac n2+s-2r\leq 0$. 
%If $r=s+3-2\alpha$, the inequality is valid for $s\geq \frac n2+4\alpha-6$.

Following similar strategy as for $I_{311}$, we estimate $I_{32}$ as follows,
\begin{equation}\notag
\begin{split}
|I_{32}|
%=&\left|\sum_{q\geq -1}\sum_{|q-p|\leq 2}\lambda_q^{2r}\int_{\R^3}\Delta_q(u_{p}\cdot\nabla b_{\leq{p-2}})b_q\, dx\right|\\
\leq& \sum_{q\geq -1}\sum_{|q-p|\leq 2}\lambda_q^{2r}\|u_p\|_2\|\nabla b_{\leq p-2}\|_\infty\|b_q\|_2\\
\lesssim & \sum_{q\geq -1}\lambda_q^{2r}\|u_q\|_2\|b_q\|_2\sum_{p\leq q}\lambda_p^{\frac n2+1}\|b_p\|_2\\
\lesssim & \sum_{q\geq -1}\lambda_q^{s+1}\|u_q\|_2\lambda_q^{(r+1)\theta}\|b_q\|_2^\theta \lambda_q^{r(1-\theta)}\|b_q\|_2^{1-\theta}\cdot\\
&\sum_{p\leq q}\lambda_p^{r}\|b_p\|_2\lambda_{q-p}^{r-s-1-\theta}\lambda_p^{\frac n2-s-\theta}\\
\lesssim & \sum_{q\geq -1}\lambda_q^{s+1}\|u_q\|_2\lambda_q^{(r+1)\theta}\|b_q\|_2^\theta \lambda_q^{r(1-\theta)}\|b_q\|_2^{1-\theta}\cdot
\sum_{p\leq q}\lambda_p^{r}\|b_p\|_2\lambda_{q-p}^{r-s-1-\theta}\\
\end{split}
\end{equation}
for $0<\theta<1$ and
\begin{equation}\label{s4}
s\geq \frac n2-\theta.
\end{equation}
It then follows from Young's inequality and Jensen's inequality, with the triplet $(2,\frac 2\theta, \frac2{1-\theta})$ satisfying
\begin{equation}\label{para7}
r-s-1-\theta<0
\end{equation}
such that
\begin{equation}\notag
\begin{split}
|I_{32}|
\lesssim & \frac{\nu}{32}\sum_{q\geq -1}\lambda_q^{2s+2}\|u_q\|_2^2+\frac{\mu}{32}\sum_{q\geq -1}\lambda_q^{2s+2}\|b_q\|_2^2
+C_{\nu,\mu}\left(\sum_{q\geq-1}\lambda_p^{2r}\|b_p\|_2^2\right)^{\frac1{1-\theta}}.
\end{split}
\end{equation}
The constraints (\ref{s4}) and (\ref{para7}) implies that for $\theta=1-\varepsilon$
\begin{equation}\notag%\label{s5}
s>r-1-\theta, \ \ s\geq \frac n2-1+\varepsilon>\frac n2-1.
\end{equation}
The term $I_{33}$ can be estimated in an analogous way as for $I_{23}$. 
To not over burden the analysis with computations, we omit the details and claim
\begin{equation}\notag
\begin{split}
|I_{33}|\leq&
\frac{\nu}{32}\sum_{q\geq -1}\lambda_q^{2s+2}\|u_q\|_2^2+\frac{\mu}{32}\sum_{q\geq -1}\lambda_q^{2s+2}\|b_q\|_2^2
+C_{\nu,\mu}\left(\sum_{q\geq-1}\lambda_p^{2r}\|b_p\|_2^2\right)^{1+\gamma_4/2}
\end{split}
\end{equation}
for some constant $\gamma_4>0$. 
%If $r=s+1$, the inequality is valid for $s\geq \frac n2-2$.
 \cbdu

 \begin{Lemma}
Let the index $r$ and $s$ satisfy conditions in Lemma \ref{le-i3}. In addition, assume $r\leq s+1-\varepsilon$ for a small enough constant $\varepsilon>0$. We have 
\begin{equation}\notag%\label{est-i4}
\begin{split}
|I_4|\leq &\frac \nu{32}\sum_{q\geq -1}\lambda_q^{2s+2}\|u_q\|_2^2+\frac \mu{32}\sum_{q\geq -1}\lambda_q^{2r+2}\|b_q\|_2^2\\
&+C_{\nu,\mu} \|u\|_{H^s}^{2+\gamma_5}
+C_{\nu,\mu} \|b\|_{H^r}^{2+\gamma_6} +C_{\nu,\mu} \|b\|_{H^r}^{2+\gamma_7}
\end{split}
\end{equation}
for various constants $C_{\nu,\mu}$ depending on $\nu,\mu$, and some constants $\gamma_5,\gamma_6,\gamma_7>0$.
 \end{Lemma}
\pf 
As usual, using Bony's paraproduct, $I_4$ can be written as
\begin{equation}\notag%\label{eq:i4}
\begin{split}
I_4=
&-\sum_{q\geq -1}\sum_{|q-p|\leq 2}\lambda_q^{2r}\int_{\R^3}\Delta_q(b_{\leq p-2}\cdot \nabla u_p) b_q\, dx\\
&-\sum_{q\geq -1}\sum_{|q-p|\leq 2}\lambda_q^{2r}\int_{\R^3}\Delta_q(b_{p}\cdot \nabla u_{\leq{p-2}}) b_q\, dx\\
&-\sum_{q\geq -1}\sum_{p\geq q-2}\lambda_q^{2r}\int_{\R^3}\Delta_q(\tilde b_p\cdot \nabla u_p) b_q\, dx\\
=&I_{41}+I_{42}+I_{43}.
\end{split}
\end{equation}
First we notice that $I_{42}$ and $I_{43}$ can be estimated as $I_{311}$ and $I_{33}$, respectively. While $I_{41}$ needs to be treated in a different way, since cancellation is not available here.  Applying H\"older's inequality and Bernstein's inequality first, we get
\begin{equation}\notag
\begin{split}
|I_{41}|\leq &\sum_{q\geq -1}\sum_{|q-p|\leq 2}\lambda_q^{2r}\|b_{\leq p-2}\|_\infty\|\nabla u_p\|_2\|b_q\|_2\\
\lesssim &\sum_{q\geq -1}\lambda_q^{2r+1}\|b_q\|_2\|u_q\|_2\sum_{p\leq q}\|b_{p}\|_\infty \\
\lesssim &\sum_{q\geq -1}\lambda_q^{2r+1}\|b_q\|_2\|u_q\|_2\sum_{p\leq q}\lambda_p^{\frac n2}\|b_{p}\|_2 \\
\lesssim & \sum_{q\geq -1}\left(\lambda_q^{(r+1)\delta}\|b_q\|_2^\delta\right)\left( \lambda_q^{r(1-\delta)}\|b_q\|_2^{1-\delta}\right)\left( \lambda_q^{(s+1)\eta}\|u_q\|_2^\eta\right)\left( \lambda_q^{s(1-\eta)}\|u_q\|_2^{1-\eta}\right)\\
&\cdot \left(\sum_{p\leq q}\lambda_p^r\|b_p\|_2\lambda_{q-p}^{r+1-s-\delta-\eta}\lambda_p^{\frac n2+1-s-\delta-\eta}\right)\\
\lesssim & \sum_{q\geq -1}\left(\lambda_q^{(r+1)\delta}\|b_q\|_2^\delta\right)\left( \lambda_q^{r(1-\delta)}\|b_q\|_2^{1-\delta}\right)\left( \lambda_q^{(s+1)\eta}\|u_q\|_2^\eta\right)\left( \lambda_q^{s(1-\eta)}\|u_q\|_2^{1-\eta}\right)\\
&\cdot \left(\sum_{p\leq q}\lambda_p^r\|b_p\|_2\lambda_{q-p}^{r+1-s-\delta-\eta}\right)
\end{split}
\end{equation}
provided that $\frac n2+1-s-\delta-\eta\leq 0$.
We apply Young's inequality with parameters $1\leq r_1,r_2,r_3,r_4,r_5\leq \infty$ satisfying 
\[\frac1{r_1}+\frac1{r_2}+\frac1{r_3}+\frac1{r_4}+\frac1{r_5}=1, \ \ r_1=\frac2\delta, \ \ r_3=\frac2\eta,\]
for some $\delta,\eta\in(0,1)$.
It yields that 
\begin{equation}\notag
\begin{split}
|I_{41}|\leq & \frac{\nu}{64}\sum_{q\geq -1}\lambda_q^{2s+2}\|u_q\|_2^2+\frac{\mu}{64}\sum_{q\geq -1}\lambda_q^{2r+2}\|b_q\|_2^2
+C_{\nu,\mu}\sum_{q\geq -1}\lambda_q^{r(1-\delta)r_2}\|b_q\|_2^{(1-\delta)r_2}\\
&+C_{\nu,\mu}\sum_{q\geq -1}\lambda_q^{s(1-\eta)r_4}\|u_q\|_2^{(1-\eta)r_4}
+C_{\nu,\mu}\sum_{q\geq -1}\left(\sum_{p\leq q}\lambda_p^r\|b_p\|_2\lambda_{q-p}^{r+1-s-\delta-\eta}\right)^{r_5}.
\end{split}
\end{equation}
Assume $r<s-1+\delta+\eta$. Using Jensen's inequality to the last term and exchanging the order of summation gives rise to
\begin{equation}\notag
\begin{split}
\sum_{q\geq -1}\left(\sum_{p\leq q}\lambda_p^r\|b_p\|_2\lambda_{q-p}^{r+1-s-\delta-\eta}\right)^{r_5}
\lesssim &\sum_{q\geq -1}\sum_{p\leq q}\lambda_p^{rr_5}\|b_p\|_2^{r_5}\lambda_{q-p}^{r+1-s-\delta-\eta}\\
\lesssim &\sum_{p\leq -1}\lambda_p^{rr_5}\|b_p\|_2^{r_5}\sum_{q\geq p}\lambda_{q-p}^{r+1-s-\delta-\eta}\\
\lesssim &\left(\sum_{p\leq -1}\lambda_p^{2r}\|b_p\|_2^2\right)^{\frac{r_5}2}.
\end{split}
\end{equation}
Thus one can choose $\delta$ and $\eta$ close enough to 1 and $r_2, r_4, r_5$ large enough such that $(1-\delta)r_2=2+\gamma_5$, $(1-\eta)r_4=2+\gamma_6$ and $r_5/2=1+\gamma_7/2$ with  $\gamma_5,\gamma_6,\gamma_7>0$. It then follows that
\begin{equation}\notag
\begin{split}
|I_{41}|\leq & \frac{\nu}{64}\sum_{q\geq -1}\lambda_q^{2s+2}\|u_q\|_2^2+\frac{\mu}{64}\sum_{q\geq -1}\lambda_q^{2r+2}\|b_q\|_2^2\\
&+C_{\nu,\mu}\|u\|_{H^s}^{2+\gamma_5}+C_{\nu,\mu}\|b\|_{H^r}^{2+\gamma_6}+C_{\nu,\mu}\|b\|_{H^r}^{2+\gamma_7}
\end{split}
\end{equation}
Indeed, one can choose $\delta+\eta=2-\varepsilon$ with $\varepsilon=\frac12[s-(\frac n2-1)]$.

\cbdu

\begin{Lemma}\label{le-i5}
Let $r>\frac n2$. Then $I_5$ satisfies 
\begin{equation}\notag%\label{est-i51}
|I_5|\leq \frac{\mu}{16}\sum_{q\geq -1}\lambda_q^{2r+2}\|b_q\|_2^2+C_{\mu}\|b\|_{H^r}^{2+\gamma_8}+C_{\mu}\|b\|_{H^r}^{2+\gamma_9}\\
\end{equation}
for  some constants $\gamma_8,\gamma_9>0$.
\end{Lemma}
\pf
Applying Bony's paraproduct first, we decompose $I_5$ to 
\begin{equation}\notag
\begin{split}
I_{5}=&\sum_{q\geq-1}\sum_{|q-p|\leq 2}\lambda_q^{2r}\int_{\mathbb R^3}\Delta_q( b_{\leq p-2}\times(\nabla\times b_p))\cdot\nabla\times b_q\, dx\\
&+\sum_{q\geq-1}\sum_{|q-p|\leq 2}\lambda_q^{2r}\int_{\mathbb R^3}\Delta_q( b_{p}\times(\nabla\times b_{\leq p-2}))\cdot\nabla\times b_q\, dx\\
&+\sum_{q\geq-1}\sum_{p\geq q-2}\lambda_q^{2r}\int_{\mathbb R^3}\Delta_q( b_{p}\times(\nabla\times \tilde b_p))\cdot\nabla\times b_q\, dx\\
=&I_{51}+I_{52}+I_{53}.
\end{split}
\end{equation}
Using the commutator notation (\ref{comm-v}), $I_{51}$ can be further decomposed as
\begin{equation}\notag
\begin{split}
I_{51}=&\sum_{q\geq-1}\sum_{|q-p|\leq 2}\lambda_q^{2r}\int_{\mathbb R^3}[\Delta_q,b_{\leq p-2}\times\nabla\times]b_p\cdot\nabla\times b_q\, dx\\
&+\sum_{q\geq-1}\lambda_q^{2r}\int_{\mathbb R^3}b_{\leq q-2}\times(\nabla\times b_q)\cdot\nabla\times b_q\, dx\\
&+\sum_{q\geq-1}\sum_{|p-q|\leq 2}\lambda_q^{2r}\int_{\mathbb R^3}(b_{\leq p-2}-b_{\leq q-2})\times(\nabla\times (b_p)_q)\cdot\nabla\times b_q\, dx\\
=&I_{511}+I_{512}+I_{513},
\end{split}
\end{equation}
where we used the fact $\sum_{q-2\leq p\leq q+2}\Delta_qb_p=b_q$.
It is clear that $I_{512}=0$ due to the cross product property. 
By the commutator estimate in Lemma  \ref{le-Hall2}, we infer 
\begin{equation}\notag
\begin{split}
|I_{511}|\lesssim &\sum_{q\geq-1}\sum_{|p-q|\leq 2}\lambda_q^{2r+1}\|\nabla b_{\leq p-2}\|_\infty\|b_p\|_2\|b_q\|_2\\
\lesssim & \sum_{q\geq -1}\lambda_q^{2r+1}\|b_q\|_2^2\sum_{p\leq q}\lambda_{p}\|b_{p}\|_\infty\\
\lesssim & \sum_{q\geq -1}\lambda_q^{2r+1}\|b_q\|_2^2\sum_{p\leq q}\lambda_{p}^{1+\frac n2}\|b_{p}\|_2\\
\lesssim & \sum_{q\geq -1}\lambda_q^{(r+1)\theta}\|b_q\|_2^{\theta}\lambda_q^{r(2-\theta)}\|b_q\|_2^{2-\theta}\sum_{p\leq q}\lambda_{p}^{(r+1)\delta}\|b_{p}\|_2^\delta \lambda_{p}^{r(1-\delta)}\|b_{p}\|_2^{1-\delta}\lambda_p^{1+\frac n2-r-\delta}\lambda_q^{1-\theta}\\
\lesssim & \sum_{q\geq -1}\lambda_q^{(r+1)\theta}\|b_q\|_2^{\theta}\lambda_q^{r(2-\theta)}\|b_q\|_2^{2-\theta}\sum_{p\leq q}\lambda_{p}^{(r+1)\delta}\|b_{p}\|_2^\delta \lambda_{p}^{r(1-\delta)}\|b_{p}\|_2^{1-\delta}\lambda_{q-p}^{1-\theta}\\
\end{split}
\end{equation}
for $0<\theta<2$, $0<\delta<1$ and 
\begin{equation}\label{para-51}
r\geq \frac n2+2-(\theta+\delta), \ \ \ 1-\theta<0.
\end{equation}
It then follows from Young's inequality with $(r_1,r_2,r_3,r_4)\in (1,\infty)^4$ satisfying
\begin{equation}\label{para-52}
\frac1{r_1}+\frac1{r_2}+\frac1{r_3}+\frac1{r_4}=1, \ \ r_1=\frac2{\theta}, \ \ r_3=\frac2{\delta}
\end{equation}
such that 
\begin{equation}\notag
\begin{split}
|I_{511}|
\leq & \frac{\mu}{16}\sum_{q\geq -1}\lambda_q^{2r+2}\|b_q\|_2^2+C_{\mu}\left(\sum_{q\geq -1}\lambda_{p}^{2r}\|b_{p}\|_2^2\right)^{1+\bar\gamma_1}+C_{\mu}\left(\sum_{q\geq -1}\lambda_{p}^{2r}\|b_{p}\|_2^2\right)^{1+\bar\gamma_2}\\
\end{split}
\end{equation}
for some constants $\bar\gamma_1,\bar\gamma_2>0$. The conditions (\ref{para-51}) and (\ref{para-52}) imply that 
\begin{equation}\label{para-53}
r\geq \frac n2+2-2+\varepsilon>\frac n2, \ \ \alpha>\frac1{\theta}=\frac1{2-\varepsilon}>\frac 12
\end{equation}
provided $\theta$ close enough to 2 and $\delta$ close enough to 0.

The term $I_{513}$ is estimated as follows,
\begin{equation}\notag
\begin{split}
|I_{513}|\leq &\sum_{q\geq-1}\sum_{|p-q|\leq 2}\lambda_q^{2r}\int_{\mathbb R^3}\left|(b_{\leq p-2}-b_{\leq q-2})\times(\nabla\times (b_p)_q)\cdot\nabla\times b_q\right|\, dx\\
\lesssim &\sum_{q\geq -1}\sum_{|p-q|\leq 2}\lambda_q^{2r}\|\nabla b_q\|_\infty\|b_{\leq p-2}-b_{\leq q-2}\|_2\|\nabla b_p\|_2\\
\lesssim &\sum_{q\geq -1}\lambda_q^{2r+\frac n2+2}\|b_q\|_2^3\\
\lesssim &\sum_{q\geq -1}\lambda_q^{(r+1)\theta}\|b_q\|_2^{\theta}\lambda_q^{r(3-\theta)}\|b_q\|_2^{3-\theta}\lambda_q^{\frac n2+2-r-\theta}\\
\lesssim &\sum_{q\geq -1}\lambda_q^{(r+1)\theta}\|b_q\|_2^{\theta}\lambda_q^{r(3-\theta)}\|b_q\|_2^{3-\theta}
\end{split}
\end{equation}
for $0<\theta<2$ and 
\begin{equation}\label{para-54}
r\geq \frac n2+2-\theta=\frac n2+2-2+\varepsilon>\frac n2
\end{equation}
provided $\theta=2-\varepsilon$ with small enough $\varepsilon$.
Thus, we have by Young's inequality that
\begin{equation}\notag
\begin{split}
|I_{513}|\leq & \frac{\mu}{16}\sum_{q\geq -1}\lambda_q^{2r+2}\|b_q\|_2^2+C_{\mu}\left(\sum_{q\geq -1}\lambda_{p}^{2r}\|b_{p}\|_2^2\right)^{1+\bar\gamma_3}
\end{split}
\end{equation}
for some constant $\bar\gamma_3>0$.

Notice that  
\begin{equation}\notag
\begin{split}
|I_{52}|=&\left|\sum_{q\geq-1}\sum_{|q-p|\leq 2}\lambda_q^{2r}\int_{\mathbb R^3}\Delta_q(\nabla\times b_{\leq p-2}\times b_{p})\cdot\nabla\times b_q\, dx\right|\\
\lesssim &\sum_{q\geq-1}\sum_{|q-p|\leq 2}\lambda_q^{2r+1}\|b_p\|_2\|\nabla b_{\leq p-2}\|_\infty\|b_q\|_2,
\end{split}
\end{equation}
thus $I_{52}$ enjoys the same estimate as for $I_{511}$.

To estimate $I_{53}$, we proceed as, by using H\"older's inequality and Bernstein's inequality  
\begin{equation}\notag
\begin{split}
|I_{53}|\leq &\sum_{q\geq -1}\sum_{p\geq q-2}\lambda_q^{2r}\int_{\mathbb R^3}|\Delta_q(b_p\times \nabla\times\tilde b_p)\cdot\nabla\times b_q|\, dx\\
\lesssim &\sum_{q\geq -1}\lambda_q^{2r}\|\nabla b_q\|_\infty\sum_{p\geq q-3}\|b_p\|_2\|\nabla b_p\|_2\\
\lesssim &\sum_{q\geq -1}\lambda_q^{2r+1+\frac{n}2}\|b_q\|_2\sum_{p\geq q-3}\lambda_p\|b_p\|_2^2\\
\lesssim &\sum_{q\geq -1}\lambda_q^{(r+1)\theta}\|b_q\|_2^{\theta}\lambda_q^{r(1-\theta)}\|b_q\|_2^{1-\theta}\sum_{p\geq q-3}\lambda_p^{(r+1)\delta}\|b_p\|_2^{\delta}\lambda_p^{r(2-\delta)}\|b_p\|_2^{2-\delta}\\
&\cdot\lambda_{p-q}^{1-2r-\delta}\lambda_q^{\frac{n}2+2-r-(\theta+\delta)}\\
\lesssim &\sum_{q\geq -1}\lambda_q^{(r+1)\theta}\|b_q\|_2^{\theta}\lambda_q^{r(1-\theta)}\|b_q\|_2^{1-\theta}\sum_{p\geq q-3}\lambda_p^{(r+1)\delta}\|b_p\|_2^{\delta}\lambda_p^{r(2-\delta)}\|b_p\|_2^{2-\delta}\lambda_{p-q}^{1-2r-\delta}\\
\end{split}
\end{equation}
for $0<\theta<1$, $0<\delta<2$ and 
\begin{equation}\label{para-55}
r\geq \frac n2+2-(\theta+\delta), \ \ 1-2r-\delta<0.
\end{equation}
Then by Young's inequality with $(r_1,r_2,r_3,r_4)\in (1,\infty)^4$ satisfying
\begin{equation}\label{para-56}
\frac1{r_1}+\frac1{r_2}+\frac1{r_3}+\frac1{r_4}=1, \ \ r_1=\frac2{\theta}, \ \ r_3=\frac2{\delta}
\end{equation}
and Jensen's inequality, we have 
\begin{equation}\notag
\begin{split}
|I_{53}|
\leq &\frac{\mu}{16}\sum_{q\geq -1}\lambda_q^{2r+2}\|b_q\|_2^2+C_{\mu}\left(\sum_{q\geq -1}\lambda_{p}^{2r}\|b_{p}\|_2^2\right)^{1+\bar\gamma_4}+C_{\mu}\left(\sum_{q\geq -1}\lambda_{p}^{2r}\|b_{p}\|_2^2\right)^{1+\bar\gamma_5}\\
\end{split}
\end{equation}
for some constants $\bar\gamma_4,\bar\gamma_5>0$. Again, (\ref{para-55}) and (\ref{para-56}) imply
\begin{equation}\notag
r>\frac n2
\end{equation}
provided $r_2, r_4$ are large enough. 
To summarize, we have for $r>\frac n2$ 
\begin{equation}\notag%\label{est-i51}
|I_5|\leq \frac{\mu}{16}\sum_{q\geq -1}\lambda_q^{2r+2}\|b_q\|_2^2+C_{\mu}\left(\sum_{q\geq -1}\lambda_{p}^{2r}\|b_{p}\|_2^2\right)^{1+\gamma_8/2}+C_{\mu}\left(\sum_{q\geq -1}\lambda_{p}^{2r}\|b_{p}\|_2^2\right)^{1+\gamma_9/2}\\
\end{equation}
for  some constants $\gamma_8,\gamma_9>0$. In fact, we can take $\gamma_8/2$ as the smallest number of $\bar\gamma_1, ...., \bar\gamma_5$ and $\gamma_9/2$ as the largest one of these constants.

\cbdu

We are ready to show the uniform estimate for $\|u(t)\|_{H^s}^2+\|b(t)\|_{H^r}^2$ on a short time interval. 

\begin{Lemma}\label{le-final}
Assume $r$ and $s$ satisfy
\[s>\frac n2-1, \ \ r>\frac n2, \ \ \frac n4+\frac s2<r\leq s+1-\varepsilon\]
for a small enough constant $\varepsilon>0$.
There exists a time $T=T(\nu,\mu,\| u_0\|_{H^{s}},\|b_0\|_{H^{r}})$ and a constant $C_{\nu,\mu}$ depending on $\nu$ and $\mu$ such that
\[\| u(t)\|_{H^{s}}^2+\|b(t)\|_{H^{r}}^2\leq C_{\nu,\mu} \left(\| u_0\|_{H^{s}}^2+\|b_0\|_{H^{r}}^2\right), \ \ \forall t\in[0,T].\]
\end{Lemma}
\pf
Combining (\ref{ineq-uq}), (\ref{ineq-bq}), and the estimates in Lemma \ref{le-i1} to Lemma \ref{le-i5}, there exist various constants $C_{\nu,\mu}$ depending on $\nu$ and $\mu$ such that
\begin{equation}\notag%\label{energy2}
\begin{split}
&\frac{d}{dt}\left(\| u\|_{H^{s}}^2+\|b\|_{H^{r}}^2\right)+\nu\|\nabla u\|_{H^{s}}^2+\mu\|\nabla b\|_{H^{r}}^2\\
%\leq &C_{\nu,\mu}\|u\|_{H^s}^{2+\gamma_1}+C_{\nu,\mu}\| u\|_{H^s}^{2+\gamma_2}+C_{\nu,\mu}\| b\|_{H^r}^{2+\gamma_3}+C_{\nu,\mu}\| b\|_{H^r}^{2+\gamma_4}\\
%&+C_\mu\|\nabla u\|_{H^{s}}\|b\|_{H^r}^2\\
\leq &C_{\nu,\mu}\left(\| u\|_{H^{s}}^2+\|b\|_{H^{r}}^2\right)^{1+\underline\gamma}+C_{\nu,\mu}\left(\| u\|_{H^{s}}^2+\|b\|_{H^{r}}^2\right)^{1+\overline\gamma}\\
%&+C_\mu\|\nabla u\|_{H^{s}}^2\left(\| u\|_{H^{s}}^2+\|b\|_{H^{r}}^2\right)\\
\end{split}
\end{equation}
with constants $\underline\gamma=\min\{\gamma_1, ..., \gamma_9\}$ and $\overline\gamma=\max\{\gamma_1, ..., \gamma_9\}$. Denote $\psi(t)=\| u(t)\|_{H^{s}}^2+\|b(t)\|_{H^{r}}^2$. Let 
\[T=\frac 12\min\left\{\frac{1}{C_{\nu,\mu}\underline\gamma \psi^{\underline\gamma}(0) }, \frac{1}{C_{\nu,\mu}\overline\gamma \psi^{\overline\gamma}(0) }\right\}.\] 
It follows from the energy inequality above that for $t\in [0,T]$, 
\begin{equation}\notag%\label{high-norm}
\begin{split}
\|u(t)\|_{H^s}^2+\|b(t)\|_{H^r}^2\leq &\frac{\|u_0\|_{H^s}^2+\|b_0\|_{H^r}^2}{\left[1-\underline\gamma C_{\nu,\mu}\left(\|u_0\|_{H^s}^2+\|b_0\|_{H^r}^2\right)^{\underline\gamma}t\right]^{1/{\underline\gamma}}}\\
&+ \frac{\|u_0\|_{H^s}^2+\|b_0\|_{H^r}^2}{\left[1-\overline\gamma C_{\nu,\mu}\left(\|u_0\|_{H^s}^2+\|b_0\|_{H^r}^2\right)^{\overline\gamma}t\right]^{1/{\overline\gamma}}}.
\end{split}
\end{equation}
It completes the proof of the lemma and concludes the proof of Theorem \ref{thm-priori}.

\bigskip

\section{Uniqueness and continuity}
\label{sec-uniq}
In this section, we establish the uniqueness of solutions stated in Theorem \ref{thm}. The continuity in time can be obtained through a rather standard procedure, see \cite{MB}; hence we omit the proof.
\begin{Theorem}
Let $\varepsilon>0$ be small enough. Assume $(u_1,b_1,p_1)$ and $(u_2,b_2,p_2)$ are solutions of (\ref{HMHD})-(\ref{initial}) in $H^s(\R^n)\times H^{s+1-\varepsilon}(\R^n)$ satisfying the estimates in Theorem \ref{thm-priori}. Then $(u_1,b_1)=(u_2,b_2)$.
\end{Theorem}
\pf The difference $(U,B,\pi)=(u_1-u_2, b_1-b_2, p_1-p_2)$ satisfies the equations
\begin{equation}\label{eq-dif}
\begin{split}
U_t+u_2\cdot \nabla U-b_2\cdot\nabla B+U\cdot\nabla u_1-B\cdot\nabla b_1+\nabla \pi=\nu\Delta U,\\
B_t+u_2\cdot \nabla B-b_2\cdot\nabla U+U\cdot\nabla b_1-B\cdot\nabla u_1
-\nabla\times((\nabla \times b_2)\times B)\\
+\nabla\times((\nabla \times B)\times b_1)=\mu\Delta B.
\end{split}
\end{equation}
The goal is to obtain a Gr\"onwall type of inequality for the $L^2$ energy of $(U,B)$.  Thus, we take inner product of the equations of $U$ and $B$ in (\ref{eq-dif}) with $U$ and $B$, respectively, to arrive at
\begin{equation}\label{energy1-dif}
\begin{split}
&\frac{d}{dt}\left(\frac12\|U\|_2^2+\frac12\|B\|_2^2\right)+\nu\|\nabla U\|_2^2+\mu\|\nabla B\|_2^2\\
=&\int_{\R^n}(b_2\cdot\nabla) B\cdot U\, dx+\int_{\R^n}(B\cdot\nabla) b_1\cdot U\, dx-\int_{\R^n}(u_2\cdot\nabla) U\cdot U\, dx\\
&-\int_{\R^n}(U\cdot\nabla) u_1\cdot U\, dx+\int_{\R^n}(b_2\cdot\nabla) U\cdot B\, dx+\int_{\R^n}(B\cdot\nabla) u_1\cdot B\, dx\\
&-\int_{\R^n}(u_2\cdot\nabla) B\cdot B\, dx-\int_{\R^n}(U\cdot\nabla) b_1\cdot B\, dx\\
&+\int_{\R^n}\nabla\times((\nabla \times b_2)\times B)\cdot B\, dx
-\int_{\R^n}\nabla\times((\nabla \times B)\times b_1)\cdot B\, dx.
\end{split}
\end{equation}
Since $(u_1,b_1)$ and  $(u_2,b_2)$ are in $H^s(\R^n)\times H^{s+1-\varepsilon}(\R^n)$ with $s>\frac n2-1$, so is $(U, B)$. Thus it can be justified that many terms on the right hand side vanish, i.e.
\begin{equation}\notag
\begin{split}
\int_{\R^n}(u_2\cdot\nabla) U\cdot U\, dx=0, \ \ \int_{\R^n}(u_2\cdot\nabla) B\cdot B\, dx=0, \\
\int_{\R^n}\nabla\times((\nabla \times B)\times b_1)\cdot B\, dx=0\\
\int_{\R^n}(b_2\cdot\nabla) B\cdot U\, dx+\int_{\R^n}(b_2\cdot\nabla) U\cdot B\, dx=0. \\
%\int_{\R^n}(B\cdot\nabla) b_1\cdot U\, dx-\int_{\R^n}(U\cdot\nabla) b_1\cdot B\, dx=0.
\end{split}
\end{equation}
We are left to estimate the five non-zero flux terms. The first one is estimated as
\begin{equation}\notag
\begin{split}
\left| \int_{\R^n}(B\cdot\nabla) b_1\cdot U\, dx\right|=&\left| \int_{\R^n}(B\cdot\nabla) U\cdot b_1\, dx\right|\\
\leq&\|B\|_2\|\nabla U\|_2\|b_1\|_\infty\\
\leq&\frac{\nu}{8}\|\nabla U\|_2^2+C_\nu\|B\|_2^2\|b_1\|_\infty^2\\
\leq&\frac{\nu}{8}\|\nabla U\|_2^2+C_\nu\|B\|_2^2\|b_1\|_{H^{s+1-\epsilon}}^2\\
\end{split}
\end{equation}
where we used the embedding $H^{s+1-\varepsilon}\subset L^\infty$ for $s+1-\varepsilon>\frac n2$ (since we can choose $\varepsilon=\frac12[s-(\frac n2-1)]$ and $s>\frac n2-1$). Analogous computation shows 
\[\left|\int_{\R^n}(U\cdot\nabla) u_1\cdot U\, dx\right| \leq\frac{\nu}{8}\|\nabla U\|_2^2+C_\nu\|U\|_2^2\|u_1\|_{H^{s+1}}^2,\]
\[\left|\int_{\R^n}(B\cdot\nabla) u_1\cdot B\, dx\right| \leq\frac{\mu}{8}\|\nabla B\|_2^2+C_\mu\|B\|_2^2\|u_1\|_{H^{s+1}}^2,\]
\[\left|\int_{\R^n}(U\cdot\nabla) b_1\cdot B\, dx\right| \leq\frac{\mu}{8}\|\nabla B\|_2^2+C_\mu\|U\|_2^2\|b_1\|_{H^{s+1-\varepsilon}}^2.\]
In the end, we estimate the Hall term as follows
\begin{equation}\notag
\begin{split}
\left|\int_{\R^n}\nabla\times((\nabla \times b_2)\times B)\cdot B\, dx\right|=&\left|\int_{\R^n}((\nabla \times b_2)\times B)\cdot \nabla\times B\, dx\right|\\
\leq &\|\nabla\times B\|_2\|\nabla\times b_2\|_\infty\|B\|_2\\
\leq &\frac{\mu}{8} \|\nabla B\|_2^2+C_\mu\|\nabla\times b_2\|_\infty^2\|B\|_2^2\\
\leq &\frac{\mu}{8} \|\nabla B\|_2^2+C_\mu\|\nabla b_2\|_{H^{s+1-\varepsilon}}^2\|B\|_2^2.
\end{split}
\end{equation}
The estimates above along with (\ref{energy1-dif}) give us 
\begin{equation}\notag
\begin{split}
&\frac{d}{dt}\left(\|U\|_2^2+\|B\|_2^2\right)+\nu\|\nabla U\|_2^2+\mu \|\nabla B\|_2^2\\
\leq &C_{\nu,\mu} \left(\|u_1\|_{H^{s+1}}^2+\|\nabla b_2\|_{H^{s+1-\varepsilon}}^2+\|b_1\|_{H^{s+1-\varepsilon}}^2\right) \left(\|U\|_2^2+\|B\|_2^2\right)\\
\leq &C_{\nu,\mu} \left(\|u_1\|_{H^{s+1}}^2+\|\nabla b_2\|_{H^{s+1-\varepsilon}}^2+C\right) \left(\|U\|_2^2+\|B\|_2^2\right).
\end{split}
\end{equation}
It follows from Gr\"onwall's inequality that
\begin{equation}\notag
\begin{split}
&\|U(t)\|_2^2+\|B(t)\|_2^2\\
\leq &\left(\|U(0)\|_2^2+\|B(0)\|_2^2\right)e^{CC_{\nu,\mu}t}\exp\left\{C_{\nu,\mu}\int_0^t\|u_1(\tau)\|_{H^{s+1}}^2+\|\nabla b_2(\tau)\|_{H^{s+1-\varepsilon}}^2\, d\tau\right\}.
\end{split}
\end{equation}
Since $U(0)=B(0)=0$, $u_1\in L^2(0,T; H^{s+1})$ and $b_2\in L^2(0,T; H^{s+2-\varepsilon})$, we infer 
\[\|U(t)\|_2^2+\|B(t)\|_2^2=0, \ \ \forall t\in[0,T].\]
\cbdu

%{\bf Acknowledgement.} The author would like to thank the anonymous referee for the careful review and valuable suggestions which have helped to improve the paper a lot. 

%\bigskip

%\Endrefs
\end{document}